\documentclass[12pt]{article}
\usepackage{amsmath, amsthm, amssymb}

\title{Congruence of
multilinear forms\footnotetext{This is the authors' version of a work that was published in Linear Algebra Appl. 418 (2006) 751--762.}}

\author{Genrich R. Belitskii%
\\ Department of Mathematics\\
Ben-Gurion University of the Negev\\
Beer-Sheva 84105, Israel\\
 genrich@cs.bgu.ac.il
\and
Vladimir V. Sergeichuk
\\ Institute of
Mathematics\\ Tereshchenkivska St.
3, Kiev, Ukraine\\
sergeich@imath.kiev.ua}
\date{}

\begin{document}

\newtheorem{theorem}{Theorem}
\newtheorem{corollary}[theorem]{Corollary}
\newtheorem{lemma}[theorem]{Lemma}

\theoremstyle{definition}
\newtheorem{definition}[theorem]{Definition}

\theoremstyle{remark}
\newtheorem{example}[theorem]{Example}
\newtheorem{remark}[theorem]{Remark}

\renewcommand{\le}{\leqslant}
\renewcommand{\ge}{\geqslant}
\newcommand{\rank}{\mathop{\rm rank}\nolimits}

\maketitle

\begin{abstract}
Let
\[
F\colon U\times\dots \times U\to
{\mathbb K},\qquad G\colon V\times\dots
\times V\to {\mathbb K}
\]
be two $n$-linear forms with $n\ge 2$
on vector spaces $U$ and $V$ over a
field ${\mathbb K}$. We say that $F$
and $G$ are symmetrically equivalent if
there exist linear bijections
$\varphi_1,\dots,\varphi_n\colon U\to
V$ such that
\[
F(u_1,\dots,u_n)=
G(\varphi_{i_1}u_1,\dots,
\varphi_{i_n}u_n)
\]
for all $u_1,\dots,u_n\in U$ and
each reordering $i_1,\dots,i_n$ of
$1,\dots,n$. The forms are said to
be congruent if
$\varphi_1=\dots=\varphi_n$.

Let $F$ and $G$ be symmetrically
equivalent. We prove that

(i) if $\mathbb K=\mathbb C$, then $F$
and $G$ are congruent;

(ii) if $\mathbb K=\mathbb R$, $
F=F_1\oplus\dots\oplus F_s\oplus 0$, $
G=G_1\oplus\dots\oplus G_r\oplus 0, $
and all summands $F_i$ and $G_j$ are
nonzero and direct-sum-indecomposable,
then $s=r$ and, after a suitable
reindexing, $F_i$ is congruent to $\pm
G_i$.


 \end{abstract}

\section{Introduction}
\label{s1}

Two matrices $A$ and $B$ over a
field $\mathbb K$ are called
\emph{congruent} if $A=S^TBS$ for
some nonsingular $S$. Two matrix
pairs $(A_1,B_1)$ and $(A_2,B_2)$
are called \emph{equivalent} if
$A_1=RA_2S$ and $B_1=RB_2S$ for some
nonsingular $R$ and $S$. Clearly, if
$A$ and $B$ are congruent, then
$(A,A^T)$ and $(B,B^T)$ are
equivalent. Quite unexpectedly, the
inverse statement holds for complex
matrices too: if $(A,A^T)$ and
$(B,B^T)$ are equivalent, then $A$
and $B$ are congruent \cite[Chapter
VI, \$\,3, Theorem 3]{mal}. This
statement was extended in
\cite{roi,ser2} to arbitrary systems
of linear mappings and bilinear
forms. In this article, we extend it
to multilinear forms.

A \emph{multilinear form} (or, more
precisely, \emph{$n$-linear form},
$n\ge 2$) on a finite dimensional
vector space $U$ over a field $\mathbb
K$ is a mapping $ F\colon U\times\dots
\times U\to {\mathbb K} $ such that
\begin{multline*}\label{0a}
F(u_1,\dots,u_{i-1},au'_i+bu''_i,
u_{i+1},\dots,u_n)
\\=aF(u_1,\dots, u'_i,
\dots,u_n) +bF(u_1,\dots,u''_i,
\dots,u_n)
\end{multline*}
for all $i\in\{1,\dots,n\}$, $a,b\in
\mathbb K$, and
$u_1,\dots,u'_i,u''_i,\dots,u_n\in
U$.

\begin{definition}
Let
\begin{equation}\label{0}
F\colon U\times\dots \times U\to
{\mathbb K}, \qquad G\colon
V\times\dots \times V\to {\mathbb K}
\end{equation}
be two $n$-linear forms.

(a) $F$ and $G$ are called
\emph{equivalent} if there exist linear
bijections
$\varphi_1,\dots,\varphi_n\colon$\!\!
$U\to V$ such that
\begin{equation*}
F(u_1,\dots,u_n)=  G(\varphi_1
u_1,\dots,\varphi_n u_n)
\end{equation*}
for all $u_1,\dots,u_n\in U$.

(b) $F$ and $G$ are called
\emph{symmetrically equivalent} if
there exist linear bijections
$\varphi_1,\dots,\varphi_n\colon U\to
V$ such that
\begin{equation}\label{1}
F(u_1,\dots,u_n)=
G(\varphi_{i_1}u_1,\dots,
  \varphi_{i_n}u_n)
\end{equation}
for all $u_1,\dots,u_n\in U$ and each
reordering $i_1,\dots,i_n$ of
$1,\dots,n$.

(c) $F$ and $G$ are called
\emph{congruent} if there exists a
linear bijection $\varphi\colon U\to V$
such that
\begin{equation*}\label{2}
F(u_1,\dots,u_n)=  G(\varphi
u_1,\dots,\varphi u_n).
\end{equation*}
for all $u_1,\dots,u_n\in U$.
\end{definition}

The \emph{direct sum} of forms
\eqref{0} is the multilinear form
\[
F\oplus G\colon (U\oplus V)\times\dots
\times (U\oplus V)\to {\mathbb K}
\]
defined as follows:
\[
(F\oplus G)(u_1+v_1,\dots,u_n+v_n):=
 F(u_1,\dots,u_n)+G(v_1,\dots,v_n)
\]
for all $u_1,\dots,u_n\in U$ and
$v_1,\dots,v_n\in V$.

We will use the internal definition:
if $F\colon U\times\dots \times U\to
{\mathbb K}$ is a multilinear form,
then $F=F_1\oplus F_2$ means that
there is a decomposition
$U=U_1\oplus U_2$ such that

(i) $F(x_1,\dots,x_n)=0$ as soon as
$x_i\in U_1$ and $x_j\in U_2$ for
some $i$ and $j$,

(ii) $F_1=F|U_1$ and $F_2=F|U_2$ are
the restrictions of $F$ to $U_1$ and
$U_2$.

A multilinear form $F\colon
U\times\dots\times U\to {\mathbb K}$ is
{\it indecomposable} if for each
decomposition $F=F_1\oplus F_2$ and the
corresponding decomposition
$U=U_1\oplus U_2$ we have $U_1=0$ or
$U_2=0$.

Our main result is the following
theorem.

\begin{theorem}\label{ther}
{\em(a)} If two multilinear forms
over $\mathbb C$ are symmetrically
equivalent, then they are congruent.

{\em(b)} If two multilinear forms $F$
and $G$ over $\mathbb R$ are
symmetrically equivalent and
\[
F=F_1\oplus\dots\oplus F_s\oplus
0,\qquad G=G_1\oplus\dots\oplus
G_r\oplus 0
\]
are their decompositions such that all
summands $F_i$ and $G_j$ are nonzero
and indecomposable, then $s=r$ and,
after a suitable reindexing, each $F_i$
is congruent to $G_i$ or $-G_i$.
\end{theorem}

The statement (a) of this theorem is
proved in the next section. We prove
(b) in the end of Section \ref{s3}
basing on Corollary \ref{colt2}, in
which we argue that every $n$-linear
form $F:U\times\dots \times U\to
{\mathbb K}$ with $n\ge 3$ over an
arbitrary field ${\mathbb K}$
decomposes into a direct sum of
indecomposable forms uniquely up to
congruence of summands. Moreover, if
$F=F_1\oplus\dots\oplus F_s\oplus 0$
is a decomposition in which
$F_1,\dots,F_s$ are nonzero and
indecomposable, and
$U=U_1\oplus\dots\oplus U_{s}\oplus
U_{0}$ is the corresponding
decomposition of $U$, then the
sequence of subspaces
$U_1+U_{0},\dots,U_s+U_{0},U_{0}$ is
determined by $F$ uniquely up to
permutations of $U_1+U_0,\dots,
U_s+U_0$.

\section{Symmetric equivalence and congruence}
\label{s2}

In this section, we prove Theorem
\ref{ther}(a) and the following
theorem, which is a weakened form of
Theorem \ref{ther}(b).

\begin{theorem} \label{t1}
If two multilinear forms $F$ and $G$
over $\mathbb R$ are symmetrically
equivalent, then there are
decompositions
\[
F=F_1\oplus F_2,\qquad G=G_1\oplus G_2
\]
such that $F_1$ is congruent to
$G_1$ and $F_2$ is congruent to
$-G_2$.
\end{theorem}

Its proof is based on two lemmas.

\begin{lemma}\label{l1}
{\em (a)} Let $T$ be a nonsingular
complex matrix having a single
eigenvalue. Then
\begin{equation*}\label{6s}
 \forall m\in {\mathbb N}\quad
 \exists f(x)\in {\mathbb C}[x]\colon
 \qquad
 f(T)^m=T^{-1}.
\end{equation*}

{\em (b)} Let $T$ be a real matrix
whose set of eigenvalues consists of
one positive real number or a pair
of distinct conjugate complex
numbers. Then
\begin{equation}\label{6a}
 \forall m\in {\mathbb N}\quad
 \exists f(x)\in {\mathbb R}[x]\colon
 \qquad
 f(T)^m=T^{-1}.
\end{equation}
\end{lemma}

\begin{proof}
(a) Let $T$ be a nonsingular complex
matrix with a single eigenvalue
$\lambda$. Since the matrix $T-\lambda
I$ is nilpotent (this follows from its
Jordan canonical form), the
substitution of $T$ for $x$ into the
Taylor expansion
\begin{multline}\label{6sss}
x^{-\frac{1}{m}}=\lambda^{-\frac{1}{m}}
+\left(-\frac{1}{m}\right)
\lambda^{-\frac{1}{m}-1}(x-\lambda)
\\+\frac{1}{2!}\left(-\frac{1}{m}\right)
\left(-\frac{1}{m}-1\right)
\lambda^{-\frac{1}{m}-2}(x-\lambda)^2+\cdots
\end{multline}
gives some matrix
\begin{equation}\label{iiuh}
f(T),\qquad f(x)\in\mathbb{C}[x],
\end{equation}
satisfying $f(T)^m=T^{-1}$.

(b) Let $T$ be a square real matrix.
If it has a single eigenvalue that
is a positive real number $\lambda$,
then all coefficients in
\eqref{6sss} are real, so the matrix
\eqref{iiuh} satisfies \eqref{6a}.

Let $T$ have only two eigenvalues
\begin{equation}\label{6aa}
\lambda=a+ib,\qquad
\bar{\lambda}=a-ib \qquad
(a,b\in\mathbb R,\ b>0).
\end{equation}
It suffices to prove \eqref{6a} for any
matrix that is similar to $T$ over
$\mathbb R$, so we may suppose that $T$
is the real Jordan matrix
\[
T=R^{-1}\begin{bmatrix}
J&0\\0&\bar{J}
\end{bmatrix}R=\begin{bmatrix}
aI+F&bI\\-bI&aI+F
\end{bmatrix},\qquad
R:=\begin{bmatrix} I&-iI\\I&iI
\end{bmatrix},
\]
in which $J=\lambda I+F$ is a direct
sum of Jordan blocks with the same
eigenvalue $\lambda$ (and so $F$ is
a nilpotent upper triangular
matrix).

It suffices to prove that
\begin{equation}\label{6b}
 \forall m\in {\mathbb N}\quad
 \exists f(x)\in {\mathbb R}[x]\colon
 \qquad
 f(J)^m=J^{-1}
\end{equation}
since such $f(x)$ satisfies \eqref{6a}:
\begin{alignat*}{2}
 f(T)^m
&=f(R^{-1}(J\oplus\bar{J})R)^m
=R^{-1}f(J\oplus\bar{J})^mR\\
&=R^{-1}(f(J)^m\oplus
\overline{f(J)}\,^m)R =R^{-1}(J\oplus
\bar{J})^{-1})R=T^{-1}.
\end{alignat*}

The matrix $F$ is nilpotent, so the
substitution of $J=\lambda I+F$ into
the Taylor expansion \eqref{6sss}
gives some matrix $g(J)$ with
$g(x)\in\mathbb{C}[x]$ satisfying
$g(J)^m=J^{-1}$. Represent $g(x)$ in
the form
\[
g(x)=g_0(x)+ig_1(x),\qquad
g_0(x),g_1(x)\in \mathbb{R}[x].
\]

It suffices to prove that $J$
reduces to $iI$ by a finite sequence
of polynomial substitutions
$$
J\longmapsto h(J),\qquad h(x)\in
{\mathbb R}[x].
$$
Indeed, their composite is some
polynomial $p(x)\in {\mathbb R}[x]$
such that $p(J)=iI$, and then
$f(x):=g_0(x)+p(x)g_1(x)\in {\mathbb
R}[x]$ satisfies \eqref{6b}:
\[
f(J)^m=\big(g_0(J)+p(J)g_1(J)\big)^m=
\big(g_0(J)+ig_1(J)\big)^m=g(J)^m=J^{-1}.
\]
 First, we replace $J$ by
$b^{-1}(J-aI)$ (see \eqref{6aa}) making
$J=iI+F$. Next, we replace $J$ by
$$\frac{3}{2}J+\frac{1}{2}J^3
=\frac{3}{2}(iI+F)
+\frac{1}{2}(-iI-3F+3iF^2+F^3)=iI+F',$$
where $F':=(3iF^2+F^3)/2$. The
degree of nilpotency of $F'$ is less
than the degree of nilpotency of
$F$; we repeat the last substitution
until obtain $iI$.
\end{proof}

\begin{definition}
Let $G\colon V\times\dots \times V\to
{\mathbb K}$ be an $n$-linear form. We
say that a linear mapping $\tau\colon
V\to V$ is $G$-{\it selfadjoint} if
\begin{equation*}\label{5}
G(v_1\dots,v_{i-1}, \tau
v_i,v_{i+1}\dots,v_n)=
G(v_1\dots,v_{j-1},\tau
v_j,v_{j+1}\dots,v_n)
\end{equation*}
for all $v_1,\dots,v_n\in V$ and all
$i$ and $j$.
\end{definition}

If $\tau$ is $G$-selfadjoint, then for
every $f(x)\in{\mathbb K}[x]$ the
linear mapping $f(\tau)$ is
$G$-selfadjoint too.

\begin{lemma}\label{l2}
Let $G\colon V\times\dots \times V\to
{\mathbb K}$ be a multilinear form over
a field ${\mathbb K}$ and let
$\tau\colon V\to V$ be a
$G$-selfadjoint linear mapping. If
\begin{equation}\label{5a}
  V=V_1\oplus\dots\oplus V_s
\end{equation}
is a decomposition of $V$ into a
direct sum of $\tau$-invariant
subspaces such that the restrictions
$\tau|V_i$ and $\tau|V_j$ of $\tau$
to $V_i$ and $V_j$ have no common
eigenvalues for all $i\ne j$, then
\begin{equation}\label{5b}
    G=G_1\oplus\dots\oplus G_s,\qquad
    G_i:=G|V_i.
\end{equation}
\end{lemma}

\begin{proof}
It suffices to consider the case $s=2$.
To simplify the formulas, we assume
that $G$ is a bilinear form. Choose
$v_1\in V_1$ and $v_2\in V_2$, we must
prove that $G(v_1,v_2)=G(v_2,v_1)=0$.

Let $f(x)$ be the minimal polynomial
of $\tau|V_2$. Since $\tau|V_1$ and
$\tau|V_2$ have no common
eigenvalues, $f(\tau|V_1)\colon
V_1\to V_1$ is a bijection, so there
exists $v_1'\in V_1$ such that
$v_1=f(\tau)v_1'$. Since $\tau$ is
$G$-selfadjoint, $f(\tau)$ is
$G$-selfadjoint too, and so
\begin{align*}
G(v_1,v_2)&=G(f(\tau)v_1',v_2)=
G(v_1',f(\tau)v_2)\\&=
G(v_1',f(\tau|V_2)v_2)=
G(v_1',0v_2)=G(v_1',0)=0.
\end{align*}
Analogously, $G(v_2,v_1)=0$.
\end{proof}

\begin{proof}[Proof of Theorem
\ref{ther}(a)]  Let $n$-linear forms
\eqref{0} over ${\mathbb K}=\mathbb
C$ be symmetrically equivalent; this
means that there exist linear
bijections
$\varphi_1,\dots,\varphi_n\colon
U\to V$ satisfying \eqref{1} for
each reordering $i_1,\dots,i_n$ of
$1,\dots,n$. Let us prove by
induction that $F$ and $G$ are
congruent. Assume that
$\varphi:=\varphi_1=\dots=\varphi_t$
for some $t<n$ and prove that there
exist linear bijections
\[
\psi_1=\dots=\psi_{t}=\psi_{t+1},
\psi_{t+2},\dots,\psi_n\colon U\to V
\]
such that
\begin{equation}\label{3}
F(u_1,\dots,u_n)=
G(\psi_{i_1}u_1,\dots,
  \psi_{i_n}u_n)
\end{equation}
for all $u_1,\dots,u_n\in U$ and
each reordering $i_1,\dots,i_n$ of
$1,\dots,n$.

By \eqref{1} and since
$\varphi_1,\dots,\varphi_n$ are
bijections, for every pair of
distinct indices $i,j$ and for all
$u_i,u_j\in U$ and
$v_1,\dots,v_{i-1},
  v_{i+1},\dots,v_{j-1},
  v_{j+1},\dots,v_n\in V$,
we have
\begin{multline}\label{4}
  G(v_1,\dots,v_{i-1},\varphi u_i,
  v_{i+1},\dots,v_{j-1},
  \varphi_{t+1}u_j,v_{j+1},\dots,v_n)
           \\=
  G(v_1,\dots,v_{i-1},\varphi_{t+1}u _i,
  v_{i+1},\dots,v_{j-1},
  \varphi u_j,v_{j+1},\dots,v_n).
\end{multline}
Denote $v_i:=\varphi_{t+1}u_i$ and
$v_j:=\varphi_{t+1}u_j$. Then
\eqref{4} takes the form
\begin{equation*}\label{4a}
  G(\dots,\varphi\varphi_{t+1}^{-1} v_i,
  \dots,v_j,\dots)=
  G(\dots,v_i,\dots,
  \varphi\varphi_{t+1}^{-1}
  v_j,\dots);
\end{equation*}
this means that the linear mapping
$
\tau:=\varphi\varphi_{t+1}^{-1}\colon
V\to V $ is $G$-selfadjoint.

Let $\lambda_1,\dots,\lambda_s$ be
all distinct eigenvalues of $\tau$
and let \eqref{5a} be the
decomposition of $V$ into the direct
sum of $\tau$-invariant subspaces
such that every $\tau_i:=\tau|V_i$
has a single eigenvalue $\lambda_i$.
Lemma \ref{l2} ensures \eqref{5b}.
For every $f_i(x)\in {\mathbb
C}[x]$, the linear mapping
$f_i(\tau_i)\colon V_i\to V_i$ is
$G_i$-selfadjoint. Using Lemma
\ref{l1}(a), we take $f_i(x)$ such
that
$f_i(\tau_i)^{t+1}=\tau_i^{-1}$.
Then
\[
\rho:=f_1(\tau_1)\oplus\dots\oplus
f_s(\tau_s)\colon  V\to V
\]
is $G$-selfadjoint and
$\rho^{t+1}=\tau^{-1}$.

Define
\begin{equation}\label{6ab}
  \psi_1=\dots=\psi_{t+1}:=\rho\varphi,
\qquad \psi_{t+2}:=\varphi_{t+2},\dots,
\psi_n:=\varphi_n.
\end{equation}
Since $\rho$ is $G$-selfadjoint and
\[
 \rho^{t+1}\varphi
=\tau^{-1}\varphi =(\varphi
\varphi_{t+1}^{-1})^{-1} \varphi=
\varphi_{t+1},
\]
we have
\begin{align*}
G(\psi_1u_1,\dots, \psi_nu_n)
 &=
G(\rho\varphi u_1,\dots, \rho\varphi
u_t, \rho\varphi u_{t+1},
\varphi_{t+2}
u_{t+2},\dots,\varphi_n u_n)
 \\&=
G(\varphi u_1,\dots, \varphi u_t,
\rho^{t+1}\varphi u_{t+1},
\varphi_{t+2}
u_{t+2},\dots,\varphi_n u_n)
 \\&=
G(\varphi_1u_1,\dots, \varphi_n
u_n)= F(u_1,\dots,u_n).
\end{align*}
So \eqref{3} holds for
$i_1=1,i_2=2,\dots,i_n=n$. The equality
\eqref{3} for an arbitrary reordering
$i_1,\dots,i_n$ of $1,\dots,n$ is
proved analogously.
\end{proof}

\begin{proof}[Proof of Theorem
\ref{t1}] Let $n$-linear forms
\eqref{0} over ${\mathbb K}=\mathbb
R$ be symmetrically equivalent; this
means that there exist linear
bijections
$\varphi_1,\dots,\varphi_n\colon
U\to V$ satisfying \eqref{1} for
each reordering $i_1,\dots,i_n$ of
$1,\dots,n$. Assume that
$\varphi:=\varphi_1=\dots=\varphi_t$
for some $t<n$. Just as in the proof
of Theorem \ref{ther}(a),
$\tau:=\varphi\varphi_{t+1}^{-1}$ is
$G$-selfadjoint. Let \eqref{5a} be
the decomposition of $V$ into the
direct sum of $\tau$-invariant
subspaces such that every
$\tau_p:=\tau|V_p$ has a single real
eigenvalue $\lambda_p$ or a pair of
conjugate complex eigenvalues
$$
\lambda_p=a_p+ib_p,\qquad
\bar{\lambda}_p=a_p-ib_p,\qquad b_p>0,
$$
and $\lambda_p\ne \lambda_q$ if
$p\ne q$. Lemma \ref{l2} ensures the
decomposition \eqref{5b}.

Define the $G$-selfadjoint linear
bijection
$$
\varepsilon=\varepsilon_11_{V_1}
\oplus\dots\oplus \varepsilon_s1_{V_s}
\colon V\to V,
$$
in which $\varepsilon_i=-1$ if
$\lambda_i$ is a negative real
number, and $\varepsilon_i=1$
otherwise. Replacing $\varphi_{t+1}$
by $\varepsilon\varphi_{t+1}$, we
obtain $\tau$ without negative real
eigenvalues. But the right-hand
member of the equality \eqref{1} may
change its sign on some subspaces
$V_p$. To preserve \eqref{1}, we
also replace $\varphi_{t+2}$ with
$\varepsilon\varphi_{t+2}$ if
$t+1<n$ and replace
$G=G_1\oplus\dots\oplus G_s$ (see
\eqref{5b}) with
\begin{equation}\label{der}
\varepsilon_1 G_1\oplus\dots\oplus
\varepsilon_s G_s
\end{equation}
if $t+1=n$.
By Lemma \ref{l1}(b), for
every $i$ there exists $f_i(x)\in
{\mathbb R}[x]$ such that
$f_i(\tau_i)^{t+1}=\tau_i^{-1}$. Define
\[
\rho=f_1(\tau_1)\oplus\dots\oplus
f_s(\tau_s)\colon  V\to V,
\]
then $\rho^{t+1}=\tau^{-1}$.
Reasoning as in the proof of Theorem
\ref{ther}(a), we find that
\eqref{3} with \eqref{der} instead
of $G$ holds for the linear mappings
\eqref{6ab}.
\end{proof}

We say that two systems of
$n$-linear forms
\[
F_1,\dots,F_s\colon U\times\dots
\times U\to {\mathbb K},\qquad
G_1,\dots,G_s\colon V\times\dots
\times V\to {\mathbb K}
\]
are \emph{equivalent} if there exist
linear bijections
$\varphi_1,\dots,\varphi_n\colon
U\to V$ such that
\begin{equation*}
F_i(u_1,\dots,u_n)= G_i(\varphi_1
u_1,\dots,\varphi_n u_n).
\end{equation*}
for each $i$ and for all
$u_1,\dots,u_n\in U$. These systems
are said to be \emph{congruent} if
$\varphi_1=\dots=\varphi_n$.

For every $n$-linear form $F$, we
construct the system of $n$-linear
forms
\begin{equation}\label{9}
{\cal S}(F)= \{F^{\sigma}\,|\,\sigma\in
S_n\},\qquad
F^{\sigma}(u_1,\dots,u_n):=
F(u_{\sigma(1)},\dots,u_{\sigma(n)}),
\end{equation}
where $S_n$ denotes the set of all
substitutions on $1,\dots,n$.

The next corollary is another form
of Theorem \ref{ther}(a).

\begin{corollary} \label{jkol}
Two multilinear forms $F$ and $G$
over $\mathbb C$ are congruent if
and only if the systems of
multilinear forms ${\cal S}(F)$ and
${\cal S}(G)$ are equivalent.
\end{corollary}

To each substitution $\sigma\in
S_n$, we assign some
$\varepsilon(\sigma)\in \{1,-1\}$.
Generalizing the notions of
symmetric and skew-symmetric
bilinear forms, we say that an
$n$-linear form $F$ is {\it
$\varepsilon$-symmetric} if
$F^{\sigma}= \varepsilon(\sigma)F$
for all $\sigma\in S_n$. If $G$ is
another $\varepsilon$-symmetric
$n$-linear form, then ${\cal S}(F)$
and ${\cal S}(G)$ are equivalent if
and only if $F$ and $G$ are
equivalent. So the next corollary
follows from Corollary \ref{jkol}.

\begin{corollary}
Two $\varepsilon$-symmetric multilinear
forms over $\mathbb C$ are equivalent
if and only if they are congruent.
\end{corollary}

\section{Direct decompositions}
\label{s3}

Every bilinear form over $\mathbb C$
or $\mathbb R$ decomposes into a
direct sum of indecomposable forms
uniquely up to congruence of
summands; see the classification of
bilinear forms in \cite{hor-ser,h-s,h-s1,ser2}. In
\cite[Theorem 2 and \S 2]{ser2} this
statement was extended to all
systems of linear mappings and
bilinear forms over $\mathbb C$ or
$\mathbb R$. The next theorem shows
that a stronger statement holds for
$n$-linear forms with $n\ge 3$ over
all fields.

\begin{theorem}\label{t2}
Let $ F\colon  U\times\dots \times
U\to {\mathbb K} $ be an $n$-linear
form with $n\ge 3$ over a field
${\mathbb K}$.

{\rm(a)} Let $F=F'\oplus 0$ and let
$F'$ have no zero direct summands.
If $U=U'\oplus U_0$ is the
corresponding decomposition of $U$,
then $U_0$ is uniquely determined by
$F$ and $F'$ is determined up to
congruence.

{\rm(b)} Let $F$ have no zero direct
summands and let
$F=F_1\oplus\dots\oplus F_s$ be its
decomposition into a direct sum of
indecomposable forms. If
$U=U_1\oplus\dots\oplus U_{s}$ is
the corresponding decomposition of
$U$, then the sequence $U_1,\dots,
U_{s}$ is determined by $F$ uniquely
up to permutations.
\end{theorem}

\begin{proof}
(a) The subspace $U_0$ is uniquely
determined by $F$ since $U_0$ is the
set of all $u\in U$ satisfying
\[
F(u,x_1,\dots,x_{n-1})=
F(x_1,u,x_2\dots,x_{n-1})=\dots=
F(x_1,\dots,x_{n-1},u)=0
\]
for all $x_1,\dots,x_{n-1}\in U$.

Let $F=F'\oplus 0=G'\oplus 0$ be two
decompositions in which $F'$ and
$G'$ have no zero direct summands,
and let $U=U'\oplus U_0=V'\oplus
U_0$ be the corresponding
decompositions of $U$. Choose bases
$u_1,\dots, u_m$ of $U'$ and
$v_1,\dots, v_m$ of $V'$ such that $
u_1-v_1,\dots, u_m-v_m$ belong to
$U_0$. Then
\begin{equation*}
F(u_{i_1},\dots,u_{i_n})=
F(v_{i_1},\dots,v_{i_n})
\end{equation*}
for all
$i_1,\dots,i_n\in\{1,\dots,m\}$, and
so the linear bijection
\[
\varphi\colon U'\longrightarrow
V',\qquad u_1\mapsto v_1\ ,\dots,\
u_{m}\mapsto v_{m},
\]
gives the congruence of $F'$ and $G'$.
\medskip

(b) Let
$
F\colon  U\times\dots \times U\to
{\mathbb K}
$
be an $n$-linear form with $n\ge 3$
that has no zero direct summands, let
\begin{equation}\label{reg}
F=F_1\oplus\dots\oplus F_s=
G_1\oplus\dots\oplus G_r
\end{equation}
be two decompositions of $F$ into
direct sums of indecomposable forms,
and let
\begin{equation}\label{regw}
U=U_1\oplus\dots\oplus U_{s}=
V_1\oplus\dots\oplus V_{r}
\end{equation}
be the corresponding decompositions
of $U$.

Put
\begin{equation}\label{loj}
d_1=\dim U_1,\ \dots,\ d_s=\dim U_s
\end{equation}
and choose two bases
\begin{equation}\label{11a}
u_1,\dots, u_m\in U_1\cup\dots\cup
U_{s},\qquad v_1,\dots, v_m\in
V_1\cup\dots\cup V_{r}
\end{equation}
of the space $U$ with the following
ordering of the first basis:
\begin{equation}\label{11b}
\text{$u_1,\dots,u_{d_1}$ is a basis of
$U_1,\quad u_{d_1+1},\dots,u_{d_1+d_2}$
is a basis of $U_2,\ \dots$}.
\end{equation}
Let $C$ be the transition matrix
from $u_1,\dots, u_m$ to $v_1,\dots,
v_m$. Partition it into $s$
horizontal and $s$ vertical strips
of sizes $d_1,d_2,$\dots,\!\!
$d_{s}$. Since $C$ is nonsingular,
by interchanging its columns (i.e.,
reindexing $v_1,\dots,v_m$) we make
nonsingular all diagonal blocks.
Changing the bases \eqref{11b}, we
make elementary transformations
within the horizontal strips of $C$
and reduce it to the form
\begin{equation}\label{12}
C =\begin{bmatrix}
I_{d_1}&C_{12}&\dots& C_{1s}\\
C_{21}&I_{d_2}&\dots& C_{2s}\\
\hdotsfor{4} \\ C_{s1}&C_{s2}&\dots&
I_{d_{s}}
\end{bmatrix}.
\end{equation}

It suffices to prove that
$u_1=v_1,\dots,u_m=v_m$, that is,
\begin{equation}\label{13}
C_{pq}=0\qquad\text{if $p\ne q$.}
\end{equation}
Indeed, by \eqref{11a} $v_1\in V_p$
for some $p$. Since $F_1$ is
indecomposable, if $d_1>1$ then
$u_1,u_2\in U_1$ and so
\begin{equation}\label{llb}
F(\dots,u_1,\dots,u_2,\dots)\ne
0\quad \text{or}\quad
F(\dots,u_2,\dots,u_1,\dots)\ne 0
\end{equation}
for some elements of $U$ denoted by
points. If \eqref{13} holds, then
$u_1=v_1$ and $u_2=v_2$. Since
$v_1\in V_p$, \eqref{llb} ensures
that $v_2\notin V_q$ for all $q\ne
p$, and so $v_2\in V_p$. This means
that $U_1\subset V_p$. Therefore,
after a suitable reindexing of
$V_1,\dots,V_s$ we obtain
$U_1\subset V_1,\dots,U_r\subset
V_r$. By \eqref{regw}, $r=s$ and
$U_1= V_1,\dots,U_r= V_r$; so the
statement (b) follows from
\eqref{llb}.

Let us prove \eqref{13}. For each
substitution $\sigma\in S_n$, the
$n$-linear form $F^{\sigma}$ defined in
\eqref{9} can be given by the
$n$-dimensional matrix

\begin{equation*}\label{13ab}
\mathbb{A}^{\sigma}
=[a^{\sigma}_{ij\dots
k}]_{i,j,\dots,k=1}^m, \qquad
a^{\sigma}_{ij\dots
k}:=F^{\sigma}(u_i,u_j,\dots,u_k),
\end{equation*}
in the basis $u_1,\dots, u_m$, or by
the $n$-dimensional matrix
$$
\mathbb{B}^{\sigma}
=[b^{\sigma}_{ij\dots
k}]_{i,j,\dots,k=1}^m, \qquad
b^{\sigma}_{ij\dots
k}:=F^{\sigma}(v_i,v_j,\dots,v_k),
$$
in the basis $v_1,\dots, v_m$. Then
for all $x_1,\dots,x_n\in U$ and
their coordinate vectors
$[x_i]=(x_{1i},\dots, x_{mi})^T$ in
the basis $u_1,\dots,u_m$, we have
\begin{equation}\label{13ac}
F^{\sigma}(x_1,\dots,x_n)=
\sum_{i,j,\dots,k=1}^m
  a^{\sigma}_{ij\dots k}
  x_{i1}x_{j2}\cdots x_{kn}.
\end{equation}
If $C=[c_{ij}]$ is the transition
matrix \eqref{12}, then
\begin{equation}\label{14}
  b^{\sigma}_{i'j'\dots k'}
  =\sum_{i,j,\dots,k=1}^m
  a^{\sigma}_{ij\dots k}
  c_{ii'}c_{jj'}\cdots c_{kk'}.
  \end{equation}

By \eqref{reg}, $a^{\sigma}_{ij\dots
k}=F^{\sigma}(u_i,u_j,\dots,u_k)\ne
0$ only if all $u_i,u_j,\dots,u_k$
belong to the same space $U_l$.
Hence $\mathbb{A}^{\sigma}$ and,
analogously, $\mathbb{B}^{\sigma}$
decompose into the direct sums of
$n$-dimensional matrices:
\begin{equation}\label{13a}
\mathbb{A}^{\sigma}=\mathbb{A}^{\sigma}_1
\oplus\dots\oplus
\mathbb{A}^{\sigma}_{s},\qquad
\mathbb{B}^{\sigma}=\mathbb{B}^{\sigma}_1
\oplus\dots\oplus
\mathbb{B}^{\sigma}_{r},
\end{equation}
in which every $\mathbb{A}^{\sigma}_i$
has size $d_i\times\dots\times d_i$ and
every $\mathbb{B}^{\sigma}_j$ has size
$\dim V_j\times\dots\times\dim V_j$.

We prove \eqref{13} using induction in
$n$.
\medskip

\noindent\emph{Base of induction:}
$n=3$. The $3$-dimensional matrices
$\mathbb{A}^{\sigma}$ and
$\mathbb{B}^{\sigma}$ can be given by
the sequences of $m$-by-$m$ matrices
\begin{gather*}\label{13bbb}
A^{\sigma}_1
=[a^{\sigma}_{ij1}]_{i,j=1}^m\ ,\
\dots,\ A^{\sigma}_m =[a^{\sigma}_{ij
m}]_{i,j=1}^m,
  \\
B^{\sigma}_1
=[b^{\sigma}_{ij1}]_{i,j=1}^m\ ,\
\dots,\ B^{\sigma}_m =[b^{\sigma}_{ij
m}]_{i,j=1}^m;
\end{gather*}
we call these matrices the {\it
layers} of $\mathbb A^{\sigma}$ and
$\mathbb{B}^{\sigma}$.
   The equality \eqref{13ac} takes the
form
\begin{equation}\label{13acj}
F^{\sigma}(x_1,x_2,x_3)
=[x_1]^T(A^{\sigma}_1x_{13}+\dots+
A^{\sigma}_mx_{m3})[x_2]
\end{equation}
for all $x_1,x_2,x_3\in U$ and their
coordinate vectors
$[x_i]=(x_{1i},\dots, x_{mi})^T$ in the
basis $u_1,\dots,u_m$. Put
\begin{equation}\label{13bm}
\begin{matrix}
H^{\sigma}_1 :=A^{\sigma}_1
c_{11}+\dots+ A^{\sigma}_m c_{m1}
 \\ \hdotsfor{1} \\
H^{\sigma}_m :=A^{\sigma}_1
c_{1m}+\dots+ A^{\sigma}_m c_{mm}
\end{matrix}
\end{equation}
By \eqref{14},
\[
  b^{\sigma}_{i'j'k'}
  =\sum_{i,j=1}^m(
  a^{\sigma}_{ij1}c_{1k'}+\dots+
  a^{\sigma}_{ijm}c_{mk'})
  c_{ii'}c_{jj'},
\]
and so
\begin{equation}\label{kok}
B^{\sigma}_1=C^TH^{\sigma}_1C\ ,\
\dots,\ B^{\sigma}_m=C^TH^{\sigma}_mC.
\end{equation}

Partition $\{1,\dots,m\}$ into the
subsets
\begin{equation}\label{ikj}
{\cal I}_1=\{1,\dots,d_1\},\quad {\cal
I}_2=\{d_1+1,\dots,d_1+d_2\},\ \ldots
\end{equation}
(see \eqref{loj}). By \eqref{13a},
if $k\in{\cal I}_q$ for some $q$,
then the $k^{\rm{th}}$ layer of
$\mathbb{A}^{\sigma}$ has the form
\begin{equation}\label{kug}
A^{\sigma}_k=0_{d_1}\oplus\dots\oplus
0_{d_{q-1}}\oplus \tilde A^{\sigma}_{k}
\oplus 0_{d_{q+1}}\oplus\dots\oplus
0_{d_{s}},
\end{equation}
in which $\tilde A^{\sigma}_{k}$ is
$d_q$-by-$d_q$. So by \eqref{13bm} and
since all diagonal blocks of the matrix
\eqref{12} are the identity matrices,
\begin{equation}\label{kugq}
H^{\sigma}_k =\sum_{i\in{\cal I}_1}
\tilde A^{\sigma}_i
c_{ik}\oplus\dots\oplus \sum_{i\in{\cal
I}_{q-1}} \tilde A^{\sigma}_i
c_{ik}\oplus \tilde
A^{\sigma}_k\oplus\sum_{i\in{\cal
I}_{q+1}} \tilde A^{\sigma}_i
c_{ik}\oplus\dots\oplus \sum_{i\in{\cal
I}_{s}} \tilde A^{\sigma}_i c_{ik}.
\end{equation}

We may suppose that
\begin{equation}\label{kugqe}
\sum_{\sigma\in S_3}\sum_{i=1}^m \rank
A^{\sigma}_i\ge \sum_{\sigma\in
S_3}\sum_{i=1}^m \rank B^{\sigma}_i;
\end{equation}
otherwise we interchange the direct
sums in \eqref{reg}. By \eqref{kug} and
\eqref{kok},
\begin{equation}\label{kugoqe}
\sum_{\sigma\in S_3}\sum_{i=1}^m \rank
\tilde A^{\sigma}_i\ge \sum_{\sigma\in
S_3}\sum_{i=1}^m \rank H^{\sigma}_i;
\end{equation}

Let us fix distinct $p$ and $q$ and
prove that $C_{pq}=0$ in \eqref{12}.
Due to \eqref{kugq}, \eqref{kugoqe},
and \eqref{kug},
\begin{equation}\label{18}
\forall k\in {\cal I}_q: \qquad
\sum_{i\in{\cal I}_p} A^{\sigma}_i
c_{ik}=0.
\end{equation}
Replacing in this sum each
$A^{\sigma}_i$ by the basis vector
$u_i$, we define
\begin{equation}\label{18km}
u:=\sum_{i\in{\cal I}_p} u_i c_{ik}\in
U_p.
\end{equation}
Since
\[
[u]=(0,\dots,0,c_{d+1,k},\dots,
c_{d+d_p,k},0,\dots,0)^T,\qquad
d:=d_1+\dots+d_{p-1},
\]
by \eqref{13acj} and \eqref{18} we
have $F^{\sigma}(x,y,u)=0$ for all
$x,y\in U_p$. This equality holds
for all substitutions $\sigma\in
S_3$, hence
\begin{equation}\label{18lj}
F(u,x,y)=F(x,u,y) =F(x,y,u)=0,
\end{equation}
and so $F|u\mathbb K$ is a zero
direct summand of $F_p=F|U_p$. Since
$F_p$ is indecomposable, $u=0$; that
is, $ c_{d+1,k}=\dots=
c_{d+d_p,k}=0. $ These equalities
hold for all $k\in {\cal I}_q$,
hence $C_{pq}=0$. This proves
\eqref{13} for $n=3$.
\medskip

\noindent\emph{Induction step.} Let
$n\ge 4$ and assume that \eqref{13}
holds for all $(n-1)$-linear forms.

The $n$-dimensional matrices
$\mathbb{A}^{\sigma}$ and
$\mathbb{B}^{\sigma}$ can be given
by the sequences of
$(n-1)$-dimensional matrices
\begin{gather*}\label{13bb}
A^{\sigma}_1 =[a^{\sigma}_{i\dots
j1}]_{i,\dots,j=1}^m\ ,\ \dots,\
A^{\sigma}_m =[a^{\sigma}_{i\dots j
m}]_{i,\dots,j=1}^m,
    \\
B^{\sigma}_1 =[b^{\sigma}_{i\dots
j1}]_{i,\dots,j=1}^m\ ,\ \dots,\
B^{\sigma}_m =[b^{\sigma}_{i\dots j
m}]_{i,\dots,j=1}^m.
\end{gather*}
By \eqref{14},
\begin{equation}\label{13bck}
\begin{matrix}
  {\displaystyle b^{\sigma}_{i'\dots j' 1}
  =\sum_{i,\dots,j}
  (a^{\sigma}_{i\dots j 1}c_{11}
  +\dots+
  a^{\sigma}_{i\dots j m}c_{m1})
  c_{ii'}\cdots c_{jj'}}
  \\ \hdotsfor{1} \\
  {\displaystyle b^{\sigma}_{i'\dots j' m}
  =\sum_{i,\dots,j}
  (a^{\sigma}_{i\dots j 1}c_{1m}
  +\dots+
  a^{\sigma}_{i\dots j m}c_{mm})
  c_{ii'}\cdots c_{jj'}}
\end{matrix}
\end{equation}
Due to \eqref{13a} and analogous to
\eqref{kug}, each $A^{\sigma}_k$
with $k\in{\cal I}_q$ (see
\eqref{ikj}) is a direct sum of
$d_1\times\dots\times
d_1,\dots,d_{s}\times\dots\times
d_{s}$ matrices, and only the
$q^{\text{th}}$ summand $\tilde
A^{\sigma}_k$ may be nonzero. This
implies \eqref{kugq} for each $k$
and for $H^{\sigma}_k$ defined in
\eqref{13bm}.

For each $(n-1)$-linear form $G$,
denote by $s(G)$ the number of
\emph{nonzero} summands in a
decomposition of $G$ into a direct
sum of indecomposable forms; this
number is uniquely determined by $G$
due to induction hypothesis. Put
$s(M):=s(G)$ if $G$ is given by an
$(n-1)$-dimensional matrix $M$. By
\eqref{13bck}, the set of
$(n-1)$-linear forms given by
$(n-1)$-dimensional matrices
\eqref{13bm} is congruent to the set
of $(n-1)$-linear forms given by
$B^{\sigma}_1,\dots,B^{\sigma}_m$.
Hence
\begin{equation}\label{13bc}
 s(H^{\sigma}_1)=s(B^{\sigma}_1)\ ,\
 \dots,\
 s(H^{\sigma}_m)=s(B^{\sigma}_m).
\end{equation}

We suppose that
\[
\sum_{\sigma\in S_n}\sum_{k=1}^m
s(A^{\sigma}_k)\ge\sum_{\sigma\in
S_n}\sum_{k=1}^m s(B^{\sigma}_k),
\]
otherwise we interchange the direct
sums in \eqref{reg}. Then by
\eqref{13bc}
\begin{equation}\label{16}
\sum_{\sigma\in S_n}\sum_{k=1}^m
s(\tilde A^{\sigma}_k) \ge
\sum_{\sigma\in S_n}\sum_{k=1}^m
s(H^{\sigma}_k).
\end{equation}

Let us fix distinct $p$ and $q$ and
prove that $C_{pq}=0$ in \eqref{12}. By
\eqref{kugq},
\begin{equation*}\label{17a}
s(H^{\sigma}_k) =s(\tilde
A^{\sigma}_k)+ \sum_{p\ne q}
s\Bigl(\sum_{i\in{\cal I}_p} \tilde
A^{\sigma}_i c_{ik}\Bigr)
\end{equation*}
for each $k\in {\cal I}_q$.
Combining it with \eqref{16}, we
have
\begin{equation*}
\sum_{i\in{\cal I}_p} A^{\sigma}_i
c_{ik}=\sum_{i\in{\cal I}_p} \tilde
A^{\sigma}_i c_{ik}= 0
\end{equation*}
for each $k\in {\cal I}_q$. Define
$u$ by \eqref{18km}. As in
\eqref{18lj}, we obtain
\[
F(u,x,\dots,y)=F(x,u,\dots,y)
=\dots=F(x,\dots,y,u)=0
\]
for all $x,\dots,y\in U_p$ and so
$F|u\mathbb K$ is a zero direct
summand of $F_p=F|U_p$. Since $F_p$
is indecomposable, $u=0$; so
$C_{pq}=0$. This proves \eqref{13}
for $n>3$.
\end{proof}

\begin{remark}
Theorem \ref{t2}(b) does not hold
for bilinear forms: for example, the
matrix of scalar product is the
identity in each orthonormal basis
of a Euclidean space. This
distinction between bilinear and
$n$-linear forms with $n\ge 3$ may
be explained by the fact that
decomposable bilinear forms are more
frequent. Let us consider forms in a
two-dimensional vector space. To
decompose a bilinear form, we must
make zero two entries in its
$2\times 2$ matrix. To decompose a
trilinear form, we must make zero
six entries in its $2\times 2\times
2$ matrix. In both the cases, these
zeros are made by transition
matrices, which have four entries.
\end{remark}

\begin{corollary}\label{colt2}
Let $F\colon  U\times\dots \times
U\to {\mathbb K}$ be an $n$-linear
form with $n\ge 3$ over a field
${\mathbb K}$. If
\begin{equation}\label{llu}
F=F_1\oplus\dots\oplus F_s\oplus 0
\end{equation}
and the summands $F_1,\dots,F_s$ are
nonzero and indecomposable, then
these summands are determined by $F$
uniquely up to congruence. Moreover,
if $U=U_1\oplus\dots\oplus
U_{s}\oplus U_0$ is the
corresponding decomposition of $U$,
then the sequence of subspaces
\begin{equation}\label{9bb}
U_1+U_0,\ \dots,\ U_s+U_0,\ U_0
\end{equation}
is determined by $F$ uniquely up to
permutations of $U_1+U_0,\dots,
U_s+U_0$.
\end{corollary}

\begin{proof}[Proof of Theorem
\ref{ther}(b)] For $n=2$ this theorem
was proved in \cite[Section 2.1]{ser2}
(and was extended to arbitrary systems
of forms and linear mappings in
\cite[Theorem 2]{ser2}). For $n\ge 3$
this theorem follows from Theorem
\ref{t1} and Corollary \ref{colt2}.
\end{proof}

\end{document}